\documentclass[draft,12pt]{article}

\usepackage[latin1]{inputenc}
\usepackage{amsmath,amssymb}
\usepackage{latexsym}
\usepackage{amsthm}

\newtheorem{theorem}{Theorem}[section]
\newtheorem{corollary}[theorem]{Corollary}
\newtheorem{lemma}[theorem]{Lemma}
\newtheorem{proposition}[theorem]{Proposition}
\newtheorem{definition}[theorem]{Definition}

\newtheorem{remark}[theorem]{Remark}

\numberwithin{equation}{section}

\def\qed{{\hfill $\square$ \bigskip}}

\def\square{{\vcenter{\vbox{\hrule height.3pt
        \hbox{\vrule width.3pt height5pt \kern5pt
           \vrule width.3pt}
        \hrule height.3pt}}}}

  \def\sF {{\cal F}}

 \def\sN {{\cal N}} 
  
  \def\sU {{\cal U}}
\def\sV {{\cal V}}

 \def\bQ {{\mathbb Q}}

 \def\bZ {{\mathbb Z}}

\def\wt{\widetilde}

\def\E{{\mathbb E}}
\def\P{{\mathbb P}}

\def\lam{{\lambda}}
\def\angel#1{{\langle #1 \rangle}}
\def\bee{\begin{equation}}
\def\bet{\begin{theorem}}
\def\bep{\begin{proposition}}
\def\bef{\begin{proof}}
\def\bel{\begin{lemma}}
\def\bec{\begin{corollary}}
\def\bed{\begin{definition}}
\def\ber{\begin{remark}}
\def\eee{\end{equation}}
\def\eet{\end{theorem}}
\def\eep{\end{proposition}}
\def\eef{\end{proof}}
\def\eel{\end{lemma}}
\def\eec{\end{corollary}}
\def\eed{\end{definition}}
\def\eer{\end{remark}}

\def\R{{\mathbb R}}
\def\E{{{\mathbb E}\,}}
\def\P{{\mathbb P}}

\def\lam{{\lambda}}

\def\al{{\alpha}}

\def\qed{{\hfill $\square$ \bigskip}}

\def\eps{\varepsilon}

\def\angel#1{{\langle#1\rangle}}

 \def\qq {\qquad}

\def\wt{\widetilde}

\def\ni{\noindent }
\def\ms{\medskip}

\def \thalf {{\textstyle {\frac12}}}

\def\square{{\vcenter{\vbox{\hrule height.3pt
        \hbox{\vrule width.3pt height5pt \kern5pt
           \vrule width.3pt}
        \hrule height.3pt}}}}

\def\tfrac#1#2{{\textstyle {\frac{#1}{#2}}}}

\def\tlint{{- \kern-0.85em \int \kern-0.2em}}  
\def\dlint{{- \kern-1.05em \int \kern-0.4em}}  

  \def\sF {{\cal F}}

 \def\sN {{\cal N}} 
  
  \def\sU {{\cal U}}
\def\sV {{\cal V}}

 \def\bQ {{\mathbb Q}}

 \def\bZ {{\mathbb Z}}

\def\nn{{\nonumber}}

\begin{document}

\title{The rate of escape of the most visited site of Brownian motion} 
\author{Richard F. Bass}

\date{\today}

\maketitle

\begin{abstract}  
\noindent {\it Abstract:} 
Let $\{L^z_t\}$ be the jointly continuous local times of a one-dimensional
Brownian motion and let $L^*_t=\sup_{z\in \mathbb R} L^z_t$. Let $V_t$ be any point
$z$ such that $L^z_t=L^*_t$, a most visited site of Brownian motion. We prove
that if $\gamma>1$, then
\[\liminf_{t\to \infty} \frac{|V_t|}{\sqrt t/(\log t)^\gamma}=\infty, \qquad \mbox{\rm a.s.}, \]
 with an analogous result for simple random walk. This proves a
conjecture of Lifshits and Shi.

\vskip.2cm
\noindent \emph{Subject Classification: Primary 60J55; Secondary 60J65, 60G50}   
\end{abstract}

\section{Introduction}\label{S:intro}

Let $S_n$ be a simple random walk, let $N^k_n=\sum_{j=0}^n 1_{(S_j=k)}$
be the number of visits by the random walk to the point $k$ by time
$n$, and let $N^*_n=\sup_{k\in \bZ} N^k_n$. Let $\sU_n=
\{k\in \bZ: N^k_n=N^*_n\}$, the set of values $k$ where $N^k_n$ takes its 
maximum, and let $U_n$ be any element of $\sU_n$. We call $\sU_n$
the set of most visited sites of the random walk at time $n$.
 This concept was introduced in \cite{Movisi}, and  was
simultaneously and independently defined by \cite{Erdos-Revesz}, who called $U_n$ a favorite
point of the random walk. In \cite{Movisi} it was proved that $U_n$ is
transient, and in fact
\bee\label{intro-E1}
\liminf_{n\to \infty} \frac{|U_n|}{\sqrt n/(\log n)^\gamma}=\infty
\eee
if $\gamma>11$ and 
\bee\label{intro-E2}
\liminf_{n\to \infty} \frac{|U_n|}{\sqrt n/(\log n)^\gamma}=0
\eee
if $\gamma<1$.
It has been of considerable interest since that time to prove that there
exists $\gamma_0$ such that \eqref{intro-E1} holds if $\gamma>\gamma_0$
and \eqref{intro-E2} holds if $\gamma<\gamma_0$ and to find the value
 of $\gamma_0$.

One can state the analogous problem for Brownian motion, and \cite{Movisi}
used Brownian motion techniques and an invariance principle for local times
to derive the results for random walk from those of Brownian motion. Let
$\{L^z_t\}$ be the jointly continuous local times of  a Brownian motion and let
$\sV_t(\omega)$ be the set of values of $z$ where the function $z\to L^z_t(\omega)$
takes its maximum. We call $\sV_t$ the set of most visited points 
or the set of favorite points of Brownian motion at 
time $t$. 
In \cite{Movisi} it was proved that
if $V_t$ is any element of $\sV_t$, then
\bee\label{intro-E3}
\liminf_{t\to \infty} \frac{|V_t|}{\sqrt t/(\log t)^\gamma}=\infty
\eee
if $\gamma>11$ and 
\bee\label{intro-E4}
\liminf_{t\to \infty} \frac{|V_t|}{\sqrt t/(\log t)^\gamma}=0
\eee
if $\gamma<1$.

The bounds in \eqref{intro-E2} and \eqref{intro-E4} have 
been improved somewhat. Lifshits and Shi \cite{Lifshits-Shi}
proved that the lim inf is 0 when $\gamma=1$ as well as when $\gamma<1$.

In \cite{Movast} the most visited sites of symmetric stable processes of
order $\al$ for $\al>1$ were studied. As a by-product of the results there,
the value of $\gamma$ in \eqref{intro-E3} was improved from 11 to 9.

In Lifshits and Shi \cite{Lifshits-Shi} it was asserted that the 
value of $\gamma$ in \eqref{intro-E1} and \eqref{intro-E3} could be any value larger than  1, or equivalently, that $\gamma_0$ exists and is equal to 1. However, as Prof.\ Shi kindly informed
us, there is a subtle but serious error in the proof; see Remark \ref{R-error}
for details.

Marcus and Rosen \cite{Marcus-Rosen} subsequently showed that 
$\gamma$ in \eqref{intro-E3} could be any value larger than 3.

In this paper we prove that the assertion of Lifshits and Shi is correct, 
that \eqref{intro-E1} and \eqref{intro-E3} hold whenever $\gamma>1$.
See Theorems \ref{main-theorem-BM} and \ref{main-theorem-rw}.
Our method relies mainly on the Ray-Knight theorems and 
a moving boundary estimate due to Novikov \cite{Novikov}.

A few words about when $\sU_n$ and $\sV_t$ consist of more than one point 
are in order. Eisenbaum
\cite{Eisenbaum-these} and Leuridan \cite{Leuridan-these} have shown that at any time $t$ there are at most
two values where $L^z_t$ takes its maximum. Toth \cite{Toth} has shown
that for $n$ sufficiently large, depending on $\omega$, there are at most
3 values of $k$ which are most visited sites for $S_n$, and more 
recently Ding and Shen \cite{ding-shen} have shown that almost surely $\sU_n$
consists of 3 distinct points infinitely often. 
It turns out that the values of the lim inf in \eqref{intro-E1}-\eqref{intro-E4}  do not depend
on which value of the most visited site is chosen.

There are many  results on the most visited sites of Brownian motion and of various other processes. See \cite{Bertoin-Marsalle}, \cite{Csaki-Revesz-Shi}, \cite{Eisenbaum-1997}, \cite{Eisenbaum-Khoshnevisan}, \cite{Hu-Shi}, \cite{Khoshnevisan-Lewis}, \cite{Leuridan-1997}, \cite{Marcus}, \cite{Revesz}, and \cite{Shi-Toth}
 for some of these.

In Section \ref{S:prelim} we state our main theorems precisely and give some preliminaries.
Section \ref{S:se} contains some estimates on  local times and squared
Bessel processes of dimension 0.
These are used in Section \ref{S:glt} to establish a lower bound
on the supremum of local time at certain random times, and 
in Section \ref{S:ft} we move from random times to fixed times to obtain our
result for Brownian motion. Finally in Section \ref{S:rw} we prove
the result for random walks.

\section{Preliminaries}\label{S:prelim}

Let $W_t$ be a one-dimensional Brownian motion
and let $\{L^z_t\}$ be a jointly
continuous version of its local times.
Let $$L^*_t=\sup_{z\in \R} L^z_t.$$ We define the collection of most
visited sites of $W$  by
$$\sV_t=\{x\in \R: L^x_t=L^*_t\}.$$
Let $V_t^s=\inf\{|x|: x\in \sV_t\}$ and $V_t^\ell=\sup\{|x|: x\in \sV_t\}.$

Our main theorem can be stated as follows.

\bet\label{main-theorem-BM}
(1) If $\gamma>1$, then 
$$\liminf_{t\to \infty} \frac{V^s_t}{\sqrt t/(\log t)^\gamma}
=\infty, \qq \mbox{\rm a.s.}$$\\
(2) If $\gamma\le 1$, 
$$\liminf_{t\to \infty} \frac{V^\ell_t}{\sqrt t/(\log t)^\gamma}
=0, \qq \mbox{\rm a.s.}$$
\eet

We have the corresponding theorem for a simple random walk $S_n$.
Let 
$$N_n^k=\sum_{j=0}^n 1_{(S_j=k)},$$
the number of times $S_j$ is equal to $k$ up to time $n$.
Let $N^*_n=\max_{k\in \bZ} N^k_n$ and let
$$\sU_t=\{k\in \bZ: N^k_n=N^*_n\}.$$
Let $U_t^s=\inf\{|x|: x\in \sN_t\}$ and $U_t^\ell=\sup\{|x|: x\in \sN_t\}.$

Our second theorem is the following.

\bet\label{main-theorem-rw}
(1) If $\gamma>1$, then 
$$\liminf_{n\to \infty} \frac{U^s_n}{\sqrt n/(\log n)^\gamma}
=\infty, \qq \mbox{\rm a.s.}$$\\
(2) If $\gamma\le 1$, 
$$\liminf_{n\to \infty} \frac{U^\ell_n}{\sqrt n/(\log n)^\gamma}
=0, \qq \mbox{\rm a.s.}$$
\eet

A  process $X_t$ is called  the square of a Bessel process of dimension 0 started
at $x\ge 0$, denoted $BES(0)^2$, if it is the unique solution to the stochastic differential equation
$$X_t=x+2\sqrt {X_t}\, dW_t,$$
where $X_t\ge 0$ a.s.\ for each $t$ and $W$ is a one-dimensional Brownian
motion with filtration $\{\sF_t\}$. 
When $X_t$ hits 0, which it does almost surely, it then stays
there forever. $X$ has a scaling property: for $r>0$ 
and $X$ is started at $x$, the process
$\frac1{r}X_t$ has the same law as the process $X_{t/r}$
started at $x/r$.
If $Y_t$ is the nonnegative square root of $X_t$ and $x>0$, then 
 $Y$ is
the unique solution to the stochastic differential equation
$$Y_t=\sqrt{x}+W_t-\frac{1}{2Y_t}\, dt.$$
See \cite{Revuz-Yor} for details.

For any process $\xi_t$ let
\bee\label{deftau}
\tau_a=\tau^\xi_a=\inf\{t>0: \xi_t=a\},
\eee
the hitting time of $a$ by the process $\xi_t$.

Let
\bee\label{ub-E101}
T_r=T(r)=\inf\{t>0: L_t^0\ge r\},
\eee
the inverse local time at 0.

The main preliminary result we need is the following version of a special case of
the Ray-Knight theorems. See \cite{Knight}, \cite{Marcus-Rosen}, and \cite{Revuz-Yor}.

\bet\label{RK-T1}
Suppose $r>0$.
The processes $\{L^z_{T_r}, z\ge 0\}$ and $\{L^{-z}_{T_r}, z\ge 0\}$
are each $BES(0)^2$ processes with time parameter $z$ started at $r$ and are independent of each other.
\eet

We also need

\bep\label{RK-P2}
Let $0<r<s$.  The processes $\{L^z_{T_s}-L^z_{T_r}, z\ge 0\}$ and
$\{L^{-z}_{T_s}-L^{-z}_{T_r}, z\ge 0\}$ are each $BES(0)^2$ processes
started at $s-r$, are independent of each other, and are independent of 
the processes $\{L^z_{T_r}, z\ge 0\}$ and $\{L^{-z}_{T_r}, z\ge 0\}$.
\eep

\bef
Since the local time at 0 of a Brownian motion increases only when the
Brownian motion is at 0, then $W_{T_r}=0$ for all $r>0$.
Proposition \ref{RK-P2} follows easily from this,  the strong Markov property applied
at time $T_r$,
and Theorem \ref{RK-T1}.
\eef

We use the letter $c$ with or without subscripts to denote finite
positive constants whose exact value is unimportant and whose value
may change from line to line.

\begin{remark}\label{R-error}{\rm
The error in \cite{Lifshits-Shi} is that inequality (2.12) of that
paper need not hold. 
Let $a>0$. Note that $\sup_{y>a\sqrt t} L^y_t$ can be
decreasing in $t$ at 
some times  because the supremum is over decreasing sets.
This can happen even when $W_t>a\sqrt t$. Similarly, $\sup_{x<a\sqrt t}
L^x_t$ can be increasing in $t$ at some times  even when $W_t>a\sqrt t$ because the supremum is over
increasing sets.
}
\end{remark}

\section{Some estimates}\label{S:se}

Define
$$I^+(t,h)=\sup_{0\le z\le h} L^z_t.$$

\bep\label{ub-P2} Let $\theta>0$. There exists a positive real number $M$
depending on $\theta$
such that
$$\limsup_{t\to \infty} \frac{\sup_{s\le t}
[I^+(s, \sqrt t/(\log t)^\theta) -L_s^0]}
{\sqrt t \log \log t/(\log t)^{\theta/2}}\le M, 
\qq \mbox{\rm a.s.}$$
\eep

\bef          Let $A_n$ be the event
$$A_n=\Big\{ \sup_{s\le 2^{n+1}}[I^+(s, 2^{(n+1)/2}/(\log 2^n)^\theta)
-L^0_s]\ge M\frac{2^{n/2}\log \log 2^n}{(\log 2^{n+1})^{\theta/2}}\Big\},$$
where $M$ is a positive real to be chosen in a moment.
By scaling, the probability of $A_n$ is the same as the probability of
$$B_n=\Big\{ \sup_{s\le 1}[I^+(s,1/(\log 2^n)^\theta)
-L^0_s]\ge M\frac{2^{-1/2}\log \log 2^n}{(\log 2^{n+1})^{\theta/2}}\Big\}.$$

Lemma 5.2 of \cite{Movisi} says that if $\delta\le 1$ and $t\ge 1$, then
$$\P(\sup_{s\le t}\  \sup_{0\le x,y\le 1, |x-y|\le \delta}
|L^y_s-L^x_s|\ge \lam)\le \frac{c_1}{\delta} e^{-\lam/c_2 \delta^{1/2} t^{1/4}}.$$
Applying this with $t=1$, $\delta=1/(\log 2^n)^\theta$, $x=0$, and $$\lam=2^{-1/2}M\log\log 2^n/
(\log 2^{n+1})^{\theta/2},$$
and recalling $\P(A_n)=\P(B_n)$,
we see that $\P(A_n)$ is summable provided we choose $M$ large enough.
By the Borel-Cantelli lemma, $\P(A_n\ \mbox{\rm i.o.})=0$.
If $2^n\le t\le 2^{n+1}$ and $t$ is large enough (depending on $\omega$), then 
\begin{align*}
\sup_{s\le t}\ [I^+(s,\sqrt t/(\log t)^\theta)
-L^0_s]
&\le \sup_{s\le 2^{n+1}} [I^+(s,2^{(n+1)/2}/(\log 2^n)^\theta)
-L^0_s]\\
&\le M\frac{2^{n/2}\log \log 2^n}{(\log 2^{n+1})^{\theta/2}}\\
&\le M \sqrt t\log \log t/(\log t)^{\theta/2}.
\end{align*}
The proposition follows.
\eef

\bep\label{gambler}
Let $X_t$ be a $BES(0)^2$ and let $\P^x$ denote the law of $X$ started at $x$. Then
$$\P^1(\tau_0<\tau_{1+a})=\frac{a}{1+a}.$$
\eep

\proof We know $\tau_0<\infty$ a.s. 
Now $X$ is a continuous martingale, hence a time change of a Brownian motion,
and thus the hitting probabilities are the same as those for a Brownian
motion.
\qed

%
%
%

The next two propositions show that in many respects a $BES(0)^2$ is similar to a Brownian motion as long as it is not too close to 0.

\bep\label{pld2}
 For $X$ a $BES(0)^2$ and $x>0$,
$$\P^x(\inf_{ s\le t} X_s< x-\lam)\le c_1e^{-c_2\lam^2/xt}.$$
\eep

\proof 
Since $X\ge 0$, there is nothing to prove unless $\lam\le x$.
By a scaling argument, it suffices to suppose $x=1$.

We start by writing
\bee\label{ld5}
\P^1(\tau^X_{1-\lam}\le t)\le \P^1(\tau_2^X\le t)
+\P^1(\tau^X_{1-\lam}\le t, \tau^X_2>t).
\eee
To estimate the terms on the right hand side of \eqref{ld5} we use Doob's
inequality. Recalling that $dX_t=2\sqrt{X_t}\, dW_t$, we
have $d\angel{X}_t=4X_t\, dt$. 

Suppose $a>0$. Then
\begin{align*}
\P^1(\tau^X_2\le t)&=\P^1(\sup_{s\le t\land \tau^X_2}X_s\ge 2)=
\P^1(\sup_{s\le t\land \tau^X_2} a(X_s-1)\ge a)\\
&\le e^{-a} \E^1 \exp(a(X_{t\land \tau^X_2}-1)).
\end{align*}
To bound the expectation,
\begin{align*}
\E^1 \exp(a&(X_{t\land \tau^X_2}-1))\\&
=\E^1\Big[ \exp(a(X_{t\land \tau^X_2}-1)-\tfrac12 a^2 \angel{X}_{t\land\tau^X_2})
\exp(\tfrac12 a^2 \angel{X}_{t\land\tau^X_2})\Big]\\
&\le \E^1 \exp(a(X_{t\land \tau^X_2}-1)-\tfrac12 a^2 \angel{X}_{t\land\tau^X_2})
e^{4a^2t}.
\end{align*}
Setting $a=1/8t$ yields
$$\P^1(\tau^X_2\le t)\le e^{-1/16t}.$$

The second term of \eqref{ld5} is slightly more complicated, but quite
similar. Let $\wt X_t$ be $X_t$ stopped at time $\tau_2^X$ and use \eqref{deftau} to define $\tau^{\wt X}_{1-\lam}$.
Suppose $a>0$ and write
\begin{align*}
\P^1(\tau^X_{1-\lam}\le t, \tau_2^X>t)&
\le \P^1(\inf_{s\le t\land \tau_{1-\lam}^{\wt X}} (\wt X_s-1)\le -\lam)\\
&=\P^1(\sup_{s\le t\land \tau^{\wt X}_{1-\lam}}(-a(\wt X_s-1))\ge a\lam)\\
&\le e^{-a\lam} \E^1 \exp(a(-(\wt X_{t\land \tau^{\wt X}_{1-\lam}}-1)))
\end{align*}
and the expectation on the last line is equal to
$$\E^1\Big[ \exp(-a(\wt X_{t\land \tau^{\wt X}_{1-\lam}}-1)-\tfrac12 a^2 \angel{\wt X}
_{t\land \tau^{\wt X}_{1-\lam}}) \exp(\tfrac12 a^2 \angel{\wt X}_{t\land \tau^{\wt X}_{1-\lam}})\Big],$$
which is bounded by $e^{4a^2t}$. Setting $a=\lam/8t$ we see the second term 
on the right of \eqref{ld5} is bounded by $e^{-\lam^2/16t}$. 

Combining the two estimates for the terms on the
right hand side of \eqref{ld5}  and recalling that we are supposing $\lam\le 1$
yields the proposition.
\qed

Another approach to the preceding  proposition is to use the  results of 
\cite{Byc}.
\ms

\bep\label{pnov} Let $R>0$, let $X_t$ be a $BES(0)^2$, and let $g$ be a non-negative absolutely 
continuous function on $[0,R]$ with $g(0)>0$.
 Let $p>1$. Then
\begin{align}
\P^1(X_t\le 1+&g(t), 0\le t\le R)\label{rnov} \\
&\le c_1e^{c_2(p)R} \Big(\frac{g(0)}{\sqrt R}\Big)^{1/p^2}
 \exp\Big(\frac{1}{2(p-1)p}\int_0^R g'(s)^2\, ds\Big)
+c_3e^{-c_4/R}.\nn
\end{align}
\eep

\proof 
By Novikov \cite{Novikov}, Theorem 6, 
\begin{align}
\P^0(W_t\le g(t),& \,0\le t\le R) \label{enov1}\\
&\le c_1\Big(\Phi_0\Big(\frac{g(0)}{\sqrt R}
\Big)\Big)^{1/p} \exp\Big(\frac1{2(p-1)}\int_0^R g'(s)^2\, ds\
\Big),\nn
\end{align}
where $W$ is a Brownian motion, $\Phi_0(x)=2\Phi(x)-1$, and 
$\Phi(x)$ is the distribution function of a standard normal random 
variable.
Note $\Phi_0(x)\le cx$ for $x\ge 0$.
 
Let $Z$ be the unique solution to 
$$dZ_t=dW_t-a(Z_t)\, dt,$$
where $a(x)=1/2x$ for $x\ge 1/2$ and $a(x)=1$ for $x<1/2$. Let
$Y_t=X_t^{1/2}$.

We start by writing
\begin{align}
\P^1(X_t\le 1&+g(t), \,0\le t\le R)
\label{enov2} \\
&\le
\P^1(X_t\le 1+g(t), \,0\le t\le R, \tau^X_{1/4}>R)+\P^1(\tau^X_{1/4}\le R).\nn
\end{align}
The second term on the right is bounded by $c_1e^{-c_2/R}$ by Proposition
\ref{pld2}. The first term on the right is equal to 
\begin{align*}
\P^1(Y_t\le (1+g(t))^{1/2},& \,0\le t \le R, \tau_{1/2}^Y>R)\\
&\le \P^1(Y_t\le 1+\tfrac12 g(t), \,0\le t\le R, \tau^Y_{1/2}>R)\\
&=\P^1(Z_t\le 1+\tfrac12 g(t), \,0\le t\le R, \tau^Z_{1/2}>R)\\
&\le \P^1(B),
\end{align*}
where $$B=\{Z_t\le 1+\tfrac12 g(t), \,0\le t\le R\}$$
and $\tau^Z_{1/2}$ is defined by \eqref{deftau};
we use the fact that $Z_t=Y_t$ for $t<\tau^Y_{1/2}$.

Let $$M_t=\exp\Big(\int_0^t a(Z_s)\, dW_s-\tfrac12\int_0^t a(Z_s)^2\, ds\Big).$$
Let $\bQ$ be defined by $d\bQ/d\P^1=M_t$ on $\sF_t$.
By the Girsanov theorem, $Z_t=W_t-\int_0^t a(Z_s)\, ds$ is a Brownian motion 
under $\bQ$.

By H\"older's inequality,
$$\P^1(B)=\E_\bQ[M_R^{-1};B]\le (\E_\bQ M_R^{-r})^{1/r}
(\bQ(B))^{1/p},$$
where $r=p/(p-1)$.
We bound the second factor by \eqref{enov1}.

It remains to bound
\begin{align*}
\E_\bQ[M_R^{-r}]&=\E_\P^1[M_R^{1-r}]\\
&=\E_\P^1\Big[\exp\Big((1-r)\int_0^R a(Z_s)\, dW_s-\tfrac{1-r}{2}\int_0^R a(Z_s)^2\, ds\Big)
\Big]\\
&=\E_\P^1\Big[\exp\Big((1-r)\int_0^R a(Z_s)\, dW_s-\tfrac{(1-r)^2}{2}\int_0^R a(Z_s)^2\, ds\Big)
\\
&\qq \times \exp\Big(\frac{(1-r)^2-(1-r)}{2}\int_0^R a(Z_s)^2\, ds\Big)\Big]\\
&\le \exp\Big(\frac{r^2-r}{2}R\Big).
\end{align*}

Combining our estimates yields the proposition.
\qed

\section{Growth of local times}\label{S:glt}

Suppose $\eps\in (0,\tfrac12)$ and $0<\delta\le \tfrac12$.
Choose $p>1$ close to 1 so that $1/p^2\ge 1-\eps$. Choose $\beta\in
(0,\tfrac12)$ small so that 
$\beta^2/4p(p-1)<\eps/2$.
 Let
\bee\label{def-u}
U_t=L_{T_1}^t-1.
\eee
Recall that here $t$ is actually the space variable for local time.
Set $$g(t)=\begin{cases} 4\delta,& t\le 16\delta^2/\beta^2;\\
\beta\sqrt t, & t>16\delta^2/\beta^2.
\end{cases}$$
Let
\bee\label{def-a}
A=\{\exists t\in [0,\delta^\eps]: U_t\ge g(t)\}.
\eee

\bep\label{u-est}
$$\P(A^c)\le c_1\delta^{1-2\eps}.$$
\eep


\proof
We estimate the right hand side of \eqref{rnov} with $R=\delta^\eps$ and
$g(0)=4\delta$. Observe that
  $g'(t)$ is zero unless $t>16\delta^2/\beta^2$, in which 
case $g'(t)=\beta/2\sqrt t$. Hence
\begin{align*}
\frac{1}{2p(p-1)}\int_0^{\delta^\eps} g'(t)^2\, dt
&\le \frac{\beta^2}{8p(p-1)}\int_{16\delta^2/\beta^2}^1\frac1{t}\,dt\\
&=\frac{\beta^2}{4p(p-1)} \log (1/\delta)+c(p,\beta),
\end{align*}
where $c(p,\beta)$ depends on $p$ and $\beta$, but not $\delta$.

Therefore 
$$\P(A^c)\le c_1 (\delta^{1-\eps/2})^{1/p^2}(1/\delta)^{\beta^2/4p(p-1)}
+c_2e^{-c_3\delta^{-\eps}}
\le c_4\delta^{1-2\eps}.$$
\qed

For $s\in [0,1]$ let
\bee\label{def-xt}
X_t^s=L^{t}_{T(1+s)}-L^{t}_{T(1)}-s.
\eee
Let 
\bee\label{def-bs}
B_s=\{\exists t\in [0,\delta^\eps]: X_t^s\le -\tfrac14 g(t)\}.
\eee

For $U$, an estimate involving  a power of $\delta$ close to 1 is the best we can expect. However the
exponential estimate we obtain in the next proposition allows us to take the
supremum over a large number of values of $s$.

\bep\label{x-est}
For $s\in [0,\delta^\eps]$
$$\P(B_s)\le c_1\log(1/\delta)e^{-c_2/\delta^{\eps}}.$$
\eep

\proof
Let $I_0=[0,16\delta^2/\beta^2]$.
Let $M$ be the smallest positive integer such that $2^M(16\delta^2/\beta^2)$
is larger than $\delta^\eps$. For $1\le m\le M$
let $$I_m=[2^{m-1}(16\delta^2/\beta^2),2^m(16\delta^2/\beta^2)].$$
For $0\le m\le M$ let
$$C_m=\{\exists t\in I_m: X_t^s\le -\tfrac14 g(t)\}.$$

By Proposition \ref{pld2},
for $1\le m\le M$,
$$\P(C_{m})\le c_1\exp\Big(-c_2\frac{2^{m-1}\delta^2}
{s2^{m}\delta^2}\Big).$$
Because $s\le \delta^\eps$, this is bounded by
$ c_1 e^{-c_2\delta^{-\eps}}$.
Similarly
$$\P(C_0)\le c_1\exp\Big(-c_2\frac{\delta^2}{s\delta^2}\Big)
\le c_3 e^{-c_4\delta^{-\eps}}.$$

Since $M\le c\log (\delta^{\eps-2})$, 
$$\P(\cup_{m=0}^M C_m)\le c_1 \log (1/\delta)e^{-c_2\delta^{-\eps}}.$$
Observing that $B_s\subset \cup_{m=0}^M C_m$ completes the proof.
\qed

\bep\label{key-prop}
There exists $c$ such that
$$\P(\exists u\in [1,1+\delta^\eps]: (L^*_{T_u}-u)\le \delta)
\le c\delta^{2-4\eps}.$$
$c$ depends on $\eps$ but not $\delta$.
\eep

\proof
Let $J=[\delta^{\eps-1}]+1$ and let
$0=s_0<s_1<\cdots < s_J=\delta^\eps$ be points of the interval $[0, \delta^\eps]$
such that $s_{j+1}-s_j\le \delta$ for all $j$.
Let 
$$D_j=\{\sup_{t\ge 0}(U_t+X^{s_j}_t)\le 2\delta\}.$$
We know $\P(D_0)\le 2\delta$ by Proposition \ref{gambler}.

Suppose $1\le j\le J$. If $\omega\in A\cap B_{s_j}^c$, then 
there exists $t\in [0,\delta^\eps]$ such that $U_t(\omega)\ge g(t)$ but
$X_t^{s_j}(\omega)\ge -\tfrac14 g(t)$.
But then $$U_t(\omega)+X_t^{s_j}(\omega)\ge g(t)-\tfrac14 g(t)
\ge 3\delta,$$
which implies $\omega\notin D_j$. Therefore $D_j\subset A^c\cup B_{s_j}$.
It follows that $$\cup_{j=1}^J D_j\subset A^c\cup (\cup_{j=1}^J B_{s_j}).$$
Using Propositions \ref{u-est} and  \ref{x-est} and the fact that $J\le c\delta^{\eps-1}$,
we then have
\begin{align*}\P(\exists j\le J: \sup_{t\ge 0} (U_t+X_t^{s_j}) 
&\le 2\delta)
\le 2\delta + c_1\delta^{1-2\eps}+c_2\delta^{\eps-1}\log(1/\delta)
e^{-c_3\delta^{-\eps}}\\
&\le c_4\delta^{1-2\eps}.
\end{align*}
If $\sup_{x\ge 0} L^x_{T(1+s_j)}-(1+s_j)
\le 2\delta$, then
$\sup_{t\ge 0} (U_t+X_t^{s_j})\le 2\delta$,
and so
\bee\label{kpe1}
\P(\exists j\le J: \sup_{x\ge 0} L^x_{T(1+s_j)}-(1+s_j)\le 2\delta)\le 
c_4 \delta^{1-2\eps}. 
\eee

Let $L^+_t=\sup_{x>0} L^x_t$ and $L^-_t=\sup_{x<0} L^x_t$. If 
$L^*_{T(1+s_j)}-(1+s_j)\le 2\delta$, then
$$L^+_{T(1+s_j)}-(1+s_j)\le 2\delta \qq \mbox{and}\qq  L^-_{T(1+s_j)}-(1+s_j)\le 2\delta.$$
By independence, symmetry, and \eqref{kpe1},
$$\P(E)\le (c_1\delta^{1-2\eps})^2=c_2\delta^{2-4\eps},$$
where
$$E=\{\exists j\le J: L^*_{T(1+s_j)}-(1+s_j)\le 2\delta\}.$$

If $u\le \delta^\eps$ and $u\in [s_j,s_{j+1}]$, then
\begin{align*}L^*_{T(1+u)}-(1+u)&\ge L^*_{T(1+s_j)}-(1+s_j)+(s_j-u)\\
&\ge L^*_{T(1+s_j)}-(1+s_j)-\delta.
\end{align*}
We conclude that on the event $E^c$
$$L^*_{T(1+u)}-(1+u)> 2\delta-\delta=\delta.$$
Therefore
$$\P(\exists u\in [0,\delta^\eps]: L_{T(1+u)}^*-(1+u)\le \delta)\le c\delta^{2-4\eps}.$$
\qed

\bet\label{trandom}
If $\gamma>1/2$, then
$$\liminf_{t\to \infty} \frac{L^*_{T_t}-t}{t/(\log t)^\gamma}=\infty,
\qq\mbox{a.s.}$$
\eet

\proof
Let $r_K=2^K$, $a>0$, and $$\delta_K=\frac{a}{(\log r_K)^\gamma}.$$
Divide $[r_K,r_{K+1}]$ into $[\delta_K^{-\eps}]+1$
equal subintervals. Each subinterval will have length less than or
equal to $\delta_K^\eps r_K$. Let
$$F_K=\{\exists t\in [r_K,r_{K+1}]: (L^*_{T_t}-t)\le \delta_K r_K\}.$$
Then by scaling, Proposition \ref{key-prop}, and our bound on
the number of subintervals,
$$\P(F_K)\le c_1 \delta_K^{-\eps}\delta_K^{2-4\eps}=c_1\delta_K^{2-5\eps}.$$
If $\gamma>\thalf$, choose $\eps$ small enough so that $(2-5\eps)\gamma>1$.
By the Borel-Cantelli lemma, $\P(F_K\mbox{ i.o.})=0.$
This implies
$$\P\Big(L^*_{T_t}-t\le \frac{at}{(\log t)^\gamma}\mbox{ i.o.}\Big)=0.$$
Since $a$ is arbitrary, the theorem follows.
\qed

\section{From random times to fixed times}\label{S:ft}

Now we derive our results for fixed times  from Theorem
\ref{trandom}.
 For values $r$ where
$T_r$ is approximately $ r^2$, the argument is straightforward, but for other values
of $r$ a different argument is necessary to avoid an extraneous power
of logarithm. 

Let
$$I(t,h)=\sup_{|z|\le h} L_t^z.$$

\bet \label{ft-C1}
Let $\gamma>1$. There exists $\rho>0$ such that
with probability one, 
$$L^*_t> I(t,\sqrt t/(\log t)^\gamma)+\frac{c\sqrt t}{(\log t)^{\rho}}$$
for all $t$ sufficiently large.
\eet

\bef
Without loss of generality assume $\gamma\le 2$.
Choose $1/2<b<\gamma/2$ and then choose
$a<\gamma$ such that $\gamma/2-a/2>b$.
Suppose
$$T_{r-}\le t\le T_r,$$
where $T_{r-}=\lim_{s\to r-} T_s$.
Then $L_t^0=r$.

\ni\emph{Case 1.} $t\le r^2(\log r)^a$.
By \cite{Kesten}, for $t$ sufficiently large (depending on $\omega$),
$$r=L_t^0\le c\sqrt{t\log \log t},$$
so $\log r\le c\log t$.
By Proposition \ref{ub-P2} and symmetry, for sufficiently large $t$ (also depending on $\omega$),
\begin{align*}
I(t, \sqrt t/(\log t)^\gamma)-L_t^0&\le c\frac{\sqrt t\log \log t}{(\log t)^{\gamma/2}}\\
&\le c\frac{r(\log r)^{a/2}\log \log r}{(\log r)^{\gamma/2}}\\
&=c\frac{r\log \log r}{(\log r)^{\gamma/2-a/2}}.
\end{align*}

 For $r$ sufficiently large, for all $s\in [r/2,r)$, 
by Theorem \ref{trandom} we have 
$$L^*_{T_s}-s\ge \frac{s}{2(\log s)^b}.$$
Letting $s$ increase up to $r$,
\begin{align*}
L_t^*-r&\ge L^*_{T_{r-}}-r\ge \frac{r}{2(\log r)^b}\\
&\ge I(t, \sqrt t/(\log t)^\gamma)-r+c\frac{r}{(\log r)^b}\\
&\ge I(t, \sqrt t/(\log t)^\gamma)-r+c\frac{\sqrt t}{(\log t)^{b+a/2}}
\end{align*}
for $t$ sufficiently large.

\ni\emph{Case 2.} $t>r^2 (\log r)^a$.
Then 
$$L_t^0=r\le c_1\frac{\sqrt t}{(\log t)^{a/2}}.$$
By this, Proposition \ref{ub-P2}, and symmetry,
there exists $K>c_1$ such that 
$$I(t,\sqrt t/(\log t)^\gamma)\le L_t^0+K\frac{\sqrt t\log\log t}{(\log t)^{\gamma/2}}\le
2K\frac{\sqrt t}{(\log t)^{a/2}}$$
for $t$ large.
By Kesten's law of the iterated logarithm (see \cite{Kesten} and also
\cite{Csaki-Foldes}),
there exists $\kappa>0$ such that
for $t$ sufficiently large,
\begin{align*}
L_t^*&\ge \kappa\sqrt t/(\log \log t)^{1/2}\\
&\ge 
3K\frac{\sqrt t}{(\log t)^{a/2}}\ge I(t,\sqrt t/(\log t)^\gamma)
+K\frac{\sqrt t}{(\log t)^{a/2}}.
\end{align*}

In either case,
\bee\label{glt-E1}
L_t^*\ge I(t, \sqrt t/(\log t)^\gamma)+c\frac{\sqrt t}{(\log t)^{b+a/2}},
\eee
and we may take $\rho=b+a/2$.
\eef

\bef[Proof of Theorem \ref{main-theorem-BM}] 
Theorem \ref{main-theorem-BM}(2) is already known; see 
\cite{Lifshits-Shi}.
For (1), let $\gamma>1$. For large enough $t$,
$$L_t^*>I(t, \sqrt t/(\log t)^\gamma),$$
which means that $L^z_t$ takes its maximum for $z$ outside the interval
$$[-\sqrt t/(\log t)^\gamma, \sqrt t/(\log t)^\gamma].$$
Theorem \ref{main-theorem-BM}(1) now follows.
\eef

\section{Random walks}\label{S:rw}

\bef[Proof of Theorem \ref{main-theorem-rw}]
(2) follows from \cite{Lifshits-Shi}, so we only consider (1).
By the invariance principle of \cite{Revesz} we can find a simple
random walk $S_n$  and a Brownian motion $W_t$  such that for each $\eps>0$,
\bee\label{rw-E11}
\sup_{k\in \bZ} |L^k_n-N^k_n|=o(n^{1/4+\eps}), \qq \mbox{\rm a.s.}
\eee

If $\gamma>1$ and $K_n=\max_{k\in \bZ, |k|\le \sqrt n/(\log n)^\gamma} N^k_n$, 
by \eqref{rw-E11}, Lemma 5.3 of \cite{Movisi}, and Theorem \ref{ft-C1},
there exists $\rho>0$ such that
\begin{align*}
N^*_n&\ge L^*_n-cn^{1/4+\eps}\\
&\ge I(n,\sqrt n/(\log n)^\gamma)+ c_1\frac{\sqrt n}{(\log n)^{\rho}}-c_2n^{1/4+\eps}\\
&\ge K_n+ c_1\frac{\sqrt n}{(\log n)^{\rho}}-2c_2n^{1/4+\eps}\\
&> K_n
\end{align*}
for $n$ sufficiently large.
We conclude the most visited site of $S_n$ must be larger in absolute
value than $\sqrt n/(\log n)^\gamma$ for $n$ large.
\eef


\medskip

\ni {\bf Richard F. Bass}\\
Department of Mathematics\\
University of Connecticut \\
Storrs, CT 06269-3009, USA\\
{\tt r.bass@uconn.edu}
\ms

\end{document}